\documentclass[10pt]{article}
\usepackage{t1enc}
\usepackage[latin1]{inputenc}
\usepackage[english]{babel}
\usepackage{amsmath,amsthm}
\usepackage{amsfonts}
\usepackage{latexsym}
\usepackage[dvips]{graphicx}
\usepackage{graphicx}
\usepackage{fancyhdr}
\newtheorem{theorem}{Theorem}[section]
\newtheorem{lemma}[theorem]{Lemma}

\newtheorem{conjecture}[theorem]{Conjecture}

\begin{document}

%\pagestyle{fancy}
%\fancyhf{}

%\chead{\small \textbf{Please cite this paper in press as: X. Zhang, J.-L Wu, J. Yan, Degree conditions for the partition of a graph into triangles and quadrilaterals, Utilitas Mathmatica 86 (2011) 341-346}}

\title{Degree conditions for the partition of a graph into
triangles and quadrilaterals\footnotetext{E-mail: sdu.zhang@yahoo.com.cn (X. Zhang), jlwu@sdu.edu.cn (J.-L. Wu).}\thanks{This work is supported by NNSF (10631070, 10871119, 60673059) and RSDP (200804220001) of China.}}
\author{Xin Zhang,~Jian-Liang Wu\thanks {Corresponding author.}~~and Jin Yan \\
{\small School of Mathematics, Shandong University, Jinan, 250100, China}}
\date{\small This paper has been published in\\\textbf{Utilitas Mathmatica 86 (2011) 341--346}}

\maketitle
\begin{abstract}\baselinestretch
For two positive integers $r$ and $s$ with $r\geq 2s-2$, if $G$ is a
graph of order $3r+4s$ such that $d(x)+d(y)\geq 4r+4s$ for every
$xy\not\in E(G)$, then $G$ independently contains $r$ triangles and
$s$ quadrilaterals, which partially prove the El-Zahar's Conjecture.


\end{abstract}

\baselinestretch
\section{Introduction}
In this paper, all graphs are finite, simple and undirected. Let $G$
be a graph. We use $V(G)$, $E(G)$, $\delta(G)$ and $\Delta (G)$ to
denote the vertex set, the edge set, the minimum degree and the
maximum degree of a graph $G$. If $uv\in E(G)$, then $u$ is said to
be the $neighbor$ of $v$. We use $N(v)$ to denote the set of
neighbors of a vertex $v$. The $degree$ $d(v)=|N(v)|$. A
$k$-$vertex$ is a vertex of degree $k$. For a subgraph(or a subset)
$H$ of $G$, we denote $N(v,H)=N(v)\cap V(H)$ and let
$d(v,H)=|N(v,H)|$. The minimum degree sum
$\sigma_{2}(G)=\min\{d(x)+d(y)|x,y\in V(G), xy\not\in E(G)\}$(When
$G$ is a complete graph, we define $\sigma_{2}(G)=\infty$). For a
subset $U$ of $V(G)$, $G[U]$ denotes the subgraph of $G$ induced by
$U$. For subsets $L$ and $M$ of $V(G)$, if $L\cap M=\emptyset$, we
say that $L$ and $M$ are $independent$, and let $E(L,M)=\{uv\in
E(G): u\in L, v\in M\}$ and $e(L,M)=|E(L,M)|$. The graph $P_{k}$ is
a path with $k$ vertices, and $C_{k}$ a cycle with $k$ vertices. We
call $C_{3}$ a triangle and $C_{4}$ a quadrilateral. We use $mQ$ to
represent $m$ copies of graph $Q$. Other notations can be found in
\cite{Bondy}.

Degree conditions which guarantee that disjoint cycles with specified length exist in a graph, especially
small cycles, are investigated in lots of paper. El-Zahar \cite{zahar} gave the following
conjecture.

\begin{conjecture} \label{conj} Let $G$ be a graph. If $|V(G)|=n_1+\cdots+n_k$ and $\delta(G)\geq\lceil
n_{1}/2\rceil+\cdots+\lceil n_{k}/2\rceil$ where $n_i\geq 3(1\leq
i\leq k)$. Then $G$ contains $k$ disjoint cycles of length
$n_1,\cdots,n_k$, respectively.
\end{conjecture}

\noindent He also proved it for $k=2$. The earlier result given by
Corr\'{a}di and Hajnal \cite{hajnal} states that every graph of
order at least $3k$ and the minimum degree at least $2k$ contains
$k$ disjoint cycles. In fact, this result just proves Conjecture $1$
when $n_{1}=\cdots=n_{k}=3$. The case $n_{1}=\cdots=n_{k}=4$ is also
called Erd\"os conjecture \cite{erdos}. Randerath $et$ $al$
\cite{randerath} proved that if a graph $G$ has order $4k$ and
$\delta(G)\geq 2k$, then $G$ contains $k-1$ disjoint quadrilaterals
and a subgraph of order 4 with at least four edges such that all the
quadrilaterals are disjoint to the subgraph. It is very close to
Erd\"os conjecture. Other corresponding results can be found in
\cite{alon} and \cite{wang}.

Here we consider the case $n_{i}\in\{3,4\}$ of Conjecture
\ref{conj}. Aigner and Brandt \cite{aiger} proved that if $G$ is a
graph such that $|V(G)|=3r+4s$ and $\delta(G)\geq 2r+\frac{8s}{3}$,
then $G$ contains $r$ triangles and $s$ quadrilaterals, all vertex
disjoint. Brandt et al \cite{brandt} proved that for two positive
integers $r$ and $s$, if $G$ is a graph of order $n\geq 3r+4s$ and
$\sigma_{2}(G)\geq n+r$, then $G$ contains $r$ triangles and $s$
cycles with length at most 4 which are vertex disjoint. Recently,
Yan \cite{yan} improved the result and proved that if $G$ is a graph
of order $n\geq 3r+4s+3$ and $\sigma_{2}(G)\geq n+r$, then $G$
contains $r$ triangles and $s$ quadrilaterals, all vertex disjoint.

In the paper , we prove that for two positive integers $r$ and $s$,
if $G$ is a graph of order $n=3r+4s$ and $\sigma_{2}(G)\geq n+r$,
then $G\supseteq rC_{3}\cup(s-1)C_{4}\cup D$, where $D$ is a graph
of order four with at least four edges; moreover, if $r\geq 2s-2$,
then $G\supseteq rC_{3}\cup sC_{4}$. The result partially proves
Conjecture \ref{conj} under the condition $n_{i}\in\{3,4\}$ and
$r\geq 2s-2$, and at the same time, the corresponding result in
\cite{randerath} as described above is generalized.

\section{Some Useful Lemmas}

\begin{lemma}{\rm \cite{johansson}}
Let $P$ and $Q$ be two disjoint paths where $P=P_{3}$. If $Q=P_{2}$ and $|E(P,Q)|\geq 3$, then $G[V(P\cup
Q)]\supseteq C_{4}$. If $Q=P_{3}$ and $|E(P,Q)|\geq 4$, then $G[V(P\cup Q)]\supseteq C_{4}$.
\end{lemma}

\begin{lemma}{\rm \cite{randerath}} \label{c4.5}
Let $C=a_{1}a_{2}a_{3}a_{4}a_{1}$ be a quadrilateral of $G$ and $u$, $v$ be two non-adjacent vertices such
that $\{u,v\}\subseteq V(G)-V(C)$. If $d(u,C)+d(v,C)\geq 5$, then $G[V(C)\cup \{u,v\}]$ contains a
quadrilateral $C^{'}$ and an edge $e$ such that $C^{'}$ and $e$ are disjoint and $e$ is incident with exactly
one of $u$ and $v$.
\end{lemma}

\begin{lemma} \label{c3.5}
Let $C=a_{1}a_{2}a_{3}a_{1}$ be a triangle of $G$ and $u$, $v$ be two non-adjacent vertices such that
$\{u,v\}\subseteq V(G)-V(C)$. If $d(u,C)+d(v,C)\geq 5$, then $G[V(C)\cup \{u,v\}]$ contains a triangle
$C^{'}$ and an edge $e$ such that $C^{'}$ and $e$ are disjoint and $e$ is incident with exactly one of $u$
and $v$.
\end{lemma}

\begin{proof}
Without loss of generality, assume that $N(u,C)=\{a_{1},a_{2},a_{3}\}$ and $N(v,C)=\{a_{1},a_{2}\}$. Then we
choose $C^{'}=va_{1}a_{2}v$ and $e=ua_{3}$.
\end{proof}

\begin{lemma}{\rm \cite{randerath}} \label{c4.9}
Let $C=a_{1}a_{2}a_{3}a_{4}a_{1}$ be a quadrilateral of $G$ and
$M_{1}$, $M_{2}$ be two paths in $G$ with order $2$. Suppose $C$,
$M_{1}$, $M_{2}$ are disjoint and $e(C,M_{1}\cup M_{2})\geq 9$. Then
$G[V(C\cup M_{1}\cup M_{2})]\supseteq C_4\cup P_{4}$.
\end{lemma}

\begin{lemma} \label{c3.9}
Let $C=a_{1}a_{2}a_{3}a_{1}$ be a triangle of $G$ and let $M_{1}$,
$M_{2}$ be two paths in $G$ with order $2$. Suppose $C$, $M_{1}$,
$M_{2}$ are disjoint and $e(C,M_{1}\cup M_{2})\geq 9$. Then
$G[V(C\cup M_{1}\cup M_{2})]\supseteq C_3\cup D$ where $|D|=4$ and
$|E(D)|\geq 4$.
\end{lemma}

\begin{proof}
Let $M_{1}=uv$ and $M_{2}=xy$. Without loss of generality, we assume
that $d(u,C)\geq d(v,C)$, $d(x,C)\geq d(y,C)$ and $d(u,C)+d(v,C)\geq
d(x,C)+d(y,C)$. Since $d(u,C)+d(v,C)+d(x,C)+d(y,C)=e(C,M_{1}\cup
M_{2})\geq 9$, we have that $d(u,C)=3$, $d(u,C)+d(v,C)\geq 5$ and
$d(x,C)+d(y,C)\geq 3$. Suppose $d(x,C)+d(y,C)=3$. Then $d(x,C)\geq
2$ and $d(u,C)=d(v,C)=3$. Without loss of generality, we assume that
$a_1,a_2\in N(x)$. Thus we can choose $C_3=uva_3u$ and
$F=G[\{x,y,a_1,a_2\}]$. Suppose $d(x,C)+d(y,C)\geq 4$. Then $x$ and
$y$ have the same neighbor in $C$, say the neighbor is $a_1$. Thus
we can choose $C_3=xya_1x$ and $D=G[\{u,v,a_2,a_3\}]$.
\end{proof}

\begin{lemma} {\rm \cite{randerath}} \label{c4.max}
Let $C$ be a quadrilateral and $P$ a path of order $4$ in $G$ such
that $C$ and $P$ are independent. $G[V(C\cup P)]$ does not contain a
quadrilateral $C'$ and a path $P'$ of order $4$ such that $C'$ and
$P'$ are independent and $e(G[V(C')])>e(G[V(C)])$. If $e(C,P)\geq
9$, then $G[V(C\cup P)]\supseteq C_4\cup D$ where $|D|=4$ and
$|E(D)|\geq 4$.
\end{lemma}

Let $F_4$ be the graph such that $|F_4|=4$, $e(F_4)= 4$ and $F_4\not\supseteq C_{4}$. In fact, $F_4$ can be
got from a claw by adding a new edge. From now on, we always write the only 3-vertex of $F_4$ as $u_{0}$, the
two 2-vertices as $u_{1}$, $u_{2}$, the only 1-vertex as $u_{3}$, and $U=\{u_{1},u_{2},u_{3}\}$.

\begin{lemma}{\rm
\cite{wang}} \label{c4.f} Let $Q$ be a quadrilateral. If $Q\cap
F_4=\emptyset$ and $e(U,Q)\geq 9$, then $G[V(Q\cup F_4)]\supseteq
2C_{4}$.
\end{lemma}

\begin{lemma} \label{c3.f}
Let $T$ be a triangle and  $T\cap F_4=\emptyset$. If  $d(u_{3},T)\geq 2$, or $e(U,T)\geq 6$ and
$d(u_{3},T)>0$, then $G[V(T\cup F_4)]\supseteq C_{3}\cup C_{4}$.
\end{lemma}

\begin{proof}
Write $T=c_{1}c_{2}c_{3}c_{1}$. If $d(u_{3},T)\geq 2$, then $G[\{u_{1},u_{2},u_{0}\}]=C_{3}$ and
$G[\{u_{3},c_{1},c_{2},c_{3}\}]\supseteq C_{4}$. So we may assume that $d(u_{3},T)=1$, without loss of
generality, we assume that $u_{3}c_{1}\in E(G)$. Then $d(u_{1},T)+d(u_{2},T)\geq e(U,T)-d(u_{3},T)\geq 5$. By
the symmetry of $u_{1}$ and $u_{2}$, we just consider the case when $d(u_{1},T)=3$. If $u_{2}c_{1}\in E(G)$,
then we can choose $C_3=u_{1}c_{2}c_{3}u_{1}$ and $C_4=u_{2}c_{1}u_{3}u_{0}u_{2}$. Otherwise $d(u_{2},T)=2$,
$u_{2}c_{2},u_{2}c_{3}\in E(G)$, and it follows that we choose $C_3=u_{2}c_{2}c_{3}u_{2}$ and
$C_4=c_{1}u_{1}u_{0}u_{3}c_{1}$.
\end{proof}

\begin{lemma} \label{c3.f.7}
Let $T$ be a triangle. If $\ T\cap F_4=\emptyset$ and $e(U,T)\geq 7$, then $G[V(T\cup F_4)]\supseteq
C_{3}\cup C_{4}$.
\end{lemma}

\begin{proof}
Since $e(U,T)\geq 7$, $e(U,T)\geq 6$ and $d(u_{3},T)>0$. By Lemma \ref{c3.f}, the lemma is true.
\end{proof}

\section{Main Results and their Proofs}
\begin{lemma} \label{lemma}
For two positive integers $r$ and $s$, if $G$ is a graph of order $n=3r+4s$ and $\sigma_{2}(G)\geq n+r$, then
$G\supseteq rC_{3}\cup (s-1)C_{4}$.
\end{lemma}

This lemma follows from Yan's result \cite{yan} described in the introduction.

\begin{theorem} \label{thb} For two positive integers $r$ and $s$, if $G$ is a graph of order
$n=3r+4s$ and $\sigma_{2}(G)\geq n+r$, then $G\supseteq
rC_{3}\cup(s-1)C_{4}\cup D$, where $D$ is a graph of order four with
at least four edges.
\end{theorem}
\begin{proof} By Lemma \ref{lemma}, $G$ independently contains $r$
triangles $T_{1},\cdots,T_{r}$ and $s-1$ quadrilaterals $Q_{1},\cdots,Q_{s-1}$. Let
$H_{T}=G[\bigcup^{r}_{i=1}V(T_{i})]$, $H_{Q}=G[\bigcup^{s-1}_{i=1}V(Q_{i})]$, $H=H_{T}\bigcup H_{Q}$ and
$D=G-V(H)$. Then $|D|=4$.
%Certainly, we can choose $H$ and $D$ such that $D$ has edges as many as possible.

First, we can choose $D$ such that $D$ contains two independent edges. Suppose that there are two vertices
$x$ and $y$ in $V(D)$ such that $xy\not\in E(G)$ and $d(x,D)=d(y,D)=0$. Then $d(x,H)+d(y,H)\geq
\sigma_{2}(G)\geq 4r+4s>4(r+s-1)$, and it follows that there exists a cycle $C\in
\{T_{1},\cdots,T_{r},Q_{1},\cdots,Q_{s-1}\}$ in $H$ such that $e(\{x,y\},C)\geq 5$.  By Lemma \ref{c4.5} and
Lemma \ref{c3.5}, we have $G[V(C)\cup \{x,y\}]\supseteq C'\cup K$ where $|C'|=|C|$ and $K$ is an edge. We
replace $C$ by $C'$ in $H$. Thus the new $H$ independently contains $r$ triangles and $s-1$ quadrilaterals
and the new $D$ satisfies $|E(D)|\geq 1$. So we assume that $|E(D)|\geq 1$. Let $uv\in E(D)$ and
$\{z,w\}=V(D)-\{u,v\}$. If $D$ does not contain two independent edges, then $zw\not\in E(D)$ and
$e(\{z,w\},uv)\leq 2$. So $d(z,H)+d(w,H)\geq \sigma_{2}(G)-e(\{z,w\},uv)\geq 4r+4s-2>4(r+s-1)$. By the
similar argument as above, we can find a new $H$ and $D$ such that $D$ contains two independent edges.

Next, we can properly choose $H$ such that $D$ contains a path of
order $4$. Let $xy$ and $zw$ are two independent edges. If
$e(xy,zw)>0$, then $D\supseteq P_{4}$. Otherwise $\Sigma_{v\in
V(D)}d(v,D)=4$ and it follows that $e(xy\cup zw,H)\geq
2\sigma_{2}(G)-e(xy\cup zw,D)\geq 8r+8s-4>8(r+s-1)$. So there exists
a cycle $C\in \{T_{1},\cdots,T_{r},Q_{1},\cdots,Q_{s-1}\}$ in $H$
such that $e(xy\cup zw,C)\geq 9$. Then by Lemma \ref{c4.9} and Lemma
\ref{c3.9}, we have $G[xy\cup zw\cup C]\supseteq C'\cup P_{4}$ where
$|C'|=|C|$. Then we replace $C$ by $C'$ in $H$ and then $D$ contains
a path of order 4.

Finally, we can properly choose $H$ such that $D$ has at least four
edges. Now we choose $T_{1},\cdots,T_{r}$ and $Q_{1},\cdots,Q_{s-1}$
such that
$M=\Sigma^{r}_{i=1}e(G[V(T_{i})])+\Sigma^{s-1}_{i=1}e(G[V(Q_{i})])$
is maximal and $D$ contains a path of order 4. Suppose $|E(D)|=3$,
that is, $D$ is a path of order 4. Then $e(D,H)\geq
2\sigma_{2}(G)-2|E(D)|\geq 8r+8s-6>8(r+s-1)$. So there exists a
cycle $C\in \{T_{1},\cdots,T_{r},Q_{1},\cdots,Q_{s-1}\}$ in $H$ such
that $e(D,C)\geq 9$. By Lemma \ref{c4.max}, Lemma \ref{c3.9} and the
maximality of $M$, we have $G[V(C\cup D)]\supseteq C^{'}\cup D^{'}$
where $|C^{'}|=|C|$, $|D^{'}|=4$ and $e(D^{'})\geq 4$. Now we
replace $C$ by $C'$ in $H$ and then $D$ contains a subgraph of order
4 with at least four edges. Hence no matter cases, $D$ has at least
four edges. We complete the proof of the theorem.
\end{proof}

\begin{theorem} \label{tha}
For two positive integers $r$ and $s$ with $r\geq 2s-2$, if $G$ is a graph of order $n=3r+4s$ and
$\sigma_{2}(G)\geq n+r$, then $G\supseteq rC_{3}\cup sC_{4}$.
\end{theorem}
\begin{proof} By Theorem \ref{thb}, $G\supseteq T_{1}\bigcup\cdots\bigcup T_{r}\bigcup Q_{1}\bigcup\cdots\bigcup
Q_{s-1} \bigcup D$ where $T_{1},\cdots,T_{r}$ are triangles, $Q_{1},\cdots,Q_{s-1}$ are quadrilaterals, and
$D$ is a subgraph of order 4 with at least four edges. Let $H_{T}=G[\bigcup^{r}_{i=1}V(T_{i})]$,
$H_{Q}=G[\bigcup^{s-1}_{i=1}V(Q_{i})]$, $H=H_{T}\bigcup H_{Q}$.

Suppose $D\not\supseteq C_{4}$. Then $D=F_4$ and $u_{1}u_{3}\not\in
E(G)$, $u_{2}u_{3}\not\in E(G)$. So we have
$d(u_{1})+d(u_{2})+2d(u_{3})\geq 2\sigma_{2}(G)\geq 8r+8s$. If there
exists a quadrilateral $Q_{i}$ in $H_{Q}$ such that $e(U,Q_{i})\geq
9$, then we have $G[V(U\cup Q_{i})]\supseteq 2C_{4}$ by Lemma
\ref{c4.f}. This implies that $G\supseteq rC_{3}\cup sC_{4}$. So we
may assume that $e(U,H_{Q})\leq 8(s-1)$. On the other hand, if there
exists a triangle $T_{i}$ in $H_{T}$ such that $e(U,T_{i})\geq 7$,
then we have $G[V(U\cup T_{i})]\supseteq C_{3}\cup C_{4}$ by Lemma
\ref{c3.f.7}. This also implies that $G\supseteq rC_{3}\cup sC_{4}$.
So we may assume that $e(U,T_{i})\leq 6$ for every triangle $T_{i}$
in $H_{T}$. Suppose there are $r_{1}$ triangles satisfying
$e(U,T_{i})\leq 5$ and $r_{2}$ triangles satisfying $e(U,T_{i})=6$.
Then $r_{1}+r_{2}=r$, $e(U,H_{T})\leq 5r_{1}+6r_{2}$. If there are
some $T_i(1\leq i\leq r)$ such that $d(u_{3},T_{i})>1$, or
$d(u_{3},T_{i})>0$ and $e(U,T_{i})=6$, then $G[V(T_i\cup
F_4)]\supseteq C_{3}\cup C_{4}$ by Lemma \ref{c3.f} and it follows
that $G\supseteq rC_{3}\cup sC_{4}$. So we can assume
$d(u_{3},H_T)\leq r_{1}$. Here
$d(u_{3})=e(U,H_Q)+e(U,H_T)-(d(u_{1})-2)-(d(u_{2})-2)+1=e(U,H_Q)+e(U,H_T)-(d(u_{1})+d(u_{2})+2d(u_3))+2d(u_3)+5$.
$d(u_{3})\geq (8r+8s)-8(s-1)-(5r_{1}+6r_{2})-5=8r-5r_{1}-6r_{2}+3$.
So we have $d(u_{3},H_{Q})=d(u_{3})-d(u_{3},H_{T})-1\geq
(8r-5r_{1}-6r_{2}+3)-(r_{1}+1)=2r+2\geq 4s-2>4(s-1)$ which
contradicts to the fact $d(u_{3},H_{Q})\leq |H_Q|= 4(s-1)$. Hence
$D\supseteq C_{4}$.
\end{proof}

\end{document}